\documentclass[number,citesort,seceqn,dvips]{arxbj}
\usepackage{mathrsfs}

\aid{0}
\volume{17}
\issue{4}
\pubyear{2011}
\firstpage{1195}
\lastpage{1216}
\doi{10.3150/10-BEJ319}

\makeatletter
\newcommand{\R}{\mathbb{R}}
\newcommand{\ep}{\varepsilon}

\newtheorem{thhm}{Theorem}[section]
\newtheorem{lem}{Lemma}[section]
\newtheorem{cor}{Corollary}[section]
\newtheorem{prop}{Proposition}[section]

\newproclaim{defn}{Definition}[section]

\newremark{rem}{Remark}[section]

\newcommand{\F}{\mathscr{F}}

\newcommand{\e}{\mathrm{e}}

\newcommand{\bfE}{\mathbb{E}}
\newcommand{\bfP}{\mathbb{P}}
\makeatother

\begin{document}
\begin{frontmatter}

\title{Approximations of fractional Brownian motion}
\runtitle{Approximations of fractional Brownian motion}

\begin{aug}
\author[1]{\fnms{Yuqiang} \snm{Li}\thanksref{1}\ead[label=e1]{yqli@stat.ecnu.edu.cn}\corref{}}
\and
\author[2]{\fnms{Hongshuai} \snm{Dai}\thanksref{2}\ead[label=e2]{math\_dsh@yahoo.com.cn}}
\runauthor{Y. Li and H. Dai}
\address[1]{School of Finance and Statistics, East China Normal
University, Shanghai 200241, China.\\ \printead{e1}}
\address[2]{School of Mathematics, Central South University,
Changsha 410075, China.\\ \printead{e2}}
\end{aug}

\received{\smonth{1} \syear{2010}}
\revised{\smonth{8} \syear{2010}}

%
\begin{abstract}
Approximations of fractional Brownian
motion using Poisson processes whose parameter sets have the
same dimensions as the approximated processes have been studied in
the literature. In this paper, a special approximation to the
one-parameter fractional Brownian motion is constructed using a
two-parameter Poisson process. The proof involves the tightness
and identification of finite-dimensional distributions.
\end{abstract}

%
\begin{keyword}
\kwd{fractional Brownian motion}
\kwd{Poisson process}
\kwd{weak convergence}
\end{keyword}

\end{frontmatter}

\section{Introduction}

The fractional Brownian motion with Hurst index $H\in(0,1)$ is a
centered Gaussian process $B^H=\{B^H(t),t\geq0\}$ with covariance
function
%
\begin{equation}\label{Def:fBm}
\bfE[B^{H}(t)B^{H}(s)]=\tfrac{1}{2}(|s|^{2H}+|t|^{2H}-|t-s|^{2H}).
\end{equation}
It follows from (\ref{Def:fBm}) that $B^{H}$ is self-similar with index
$H$ and has stationary increments. Unless $H=1/2$ (i.e., $B^H$ is
Brownian motion), $B^{H}$ is not Markovian. Moreover, it is known that
$B^H$ has long-range dependence if $H \in(1/2, 1)$ and short-range
dependence if $H \in(0,  1/2)$ (see Samorodnitsky and Taqqu
\cite{ST94}). These properties have made $B^{H}$ not only important
theoretically, but also very popular as stochastic models in many areas
including telecommunications, biology, hydrology and finance.

Weak convergence to fractional Brownian motion has been studied
extensively since the works of Davydov \cite{Davydov70} and Taqqu
\cite{Taqqu75}. In recent years many new results on approximations of
fractional Brownian motion have been established. For example, Enriquez
\cite{E2004} showed that fractional Brownian motion can be approximated
in law by appropriately normalized correlated random walks. Meyer,
Sellan and Taqqu \cite{MST} proved that the law of $B^H$ can be
approximated by those of a random wavelet series. By extending Stroock
\cite{S82}, Bardina \textit{et al.}~\cite{BJ2003} and Delgado and Jolis
\cite{b2} have established approximations in law to fractional Brownian
motions by processes constructed using Poisson processes.

Let $\{N(t), t \ge0\}$ be a standard Poisson process, and for all $\ep>
0$, define the processes $X_\ep= \{X_\ep(t ), t \ge0\}$ by
\[
X_{\varepsilon}(t)=\frac{1}{\varepsilon}\int_{0}^{t}(-1)^{N(r/\varepsilon^{2})}\,\mathrm{d}r,\qquad
t\geq0.
\]
Stroock \cite{S82} proved that as $\ep$ tends to zero, the laws of
$X_{\ep}$ converge weakly in the Banach space $\mathcal{C}[0,1]$ (i.e.,
the space of continuous functions on $[0, 1]$) to the law of Brownian
motion. Delgado and Jolis \cite{b2} proved that every Gaussian process
of the form
\[
X_t=\int_0^1 K(t, s)\,\mathrm{d} B_s,
\]
where $B$ is a one-dimensional Brownian motion and $K$ a sufficiently
regular deterministic kernel, can be weakly approximated by the family
of processes
\[
Y^\varepsilon(t)=\frac{1}{\varepsilon}\int^1_0 K(t,
r)(-1)^{N(r/\varepsilon^2)} \,\mathrm{d}r.
\]
Their result can be applied to fractional Brownian motion. In
addition, Bardina and Jolis \cite{BJ001} proved that as $\ep$ tends to
$0$, the family of two-parameter random fields $Y_\ep$ defined by
%
\begin{equation}\label{intr-1}
Y_{\ep}(s,t)=\int_{0}^{t}\!\!\!\int_{0}^{s}\frac{1}{\ep^{2}}
\sqrt{xy}(-1)^{N(x/\ep,y/\ep)}\,\mathrm{d}x\,\mathrm{d}y
\end{equation}
converges in law in the space of continuous functions on $[0,1]^{2}$ to
the standard Brownian sheet. Bardina, Jolis and Tudor showed in
\cite{BJ2003} that as $\ep$ tends to $0$, the family of two-parameter
random fields
%
\begin{equation}\label{intr-11}
\hat{Y}_{\ep}(s,t)=\int_{0}^{1}\!\!\!\int_{0}^{1}\frac{1}{\ep^{2}} f(s, t, x,
y)\sqrt{xy}(-1)^{N(x/\ep,y/\ep)}\,\mathrm{d}x\,\mathrm{d}y
\end{equation}
converges in law to the two-parameter Gaussian process
\[
W(s, t)=\int_0^1\!\!\!\int_0^1 f(s,t, x, y)B(\mathrm{d}x, \mathrm{d}y),
\]
where $B$ is a standard
Brownian sheet and the deterministic kernel
%
\begin{equation}\label{intr-12}
f(s, t, x, y)=f_1(s, x)f_2(t, y)
\end{equation}
can be separated by the integration variables and satisfies certain
conditions. As examples, the authors include the fractional Brownian
sheet, among others. For more information, see Bardina, Jolis and
Rovira \cite{BJ00}, where an approximation to the $d$-parameter Wiener
process by a $d$-parameter Poisson process was provided, and Bardina
and Bascompte \cite{BB09}, where two independent Gaussian processes
were constructed using a unique Poisson process.

We note that in the serial works \cite{BB09,BJ001,BJ00,BJ2003,b2}, the
dimension of the parameter set is always the same for the approximating
and the approximated processes. Naturally, we will be interested in the
problem of whether we can approximate the $d$-parameter fractional
Brownian motions by $r$-parameter Poisson processes if $d\not=r$. The
purpose of this paper is to study this problem in the case of $d=1$ and
$r=2$. We find that for a special deterministic kernel function which
cannot be separated with respect to the integration variables, the
answer is affirmative. Below, we introduce the deterministic kernel
function.

In order to study a non-Gaussian and non-stable process arising as the
limit of sums of rescaled renewal processes under the condition of
intermediate growth, Gaigalas \cite{G06}, page 451, introduced the
function $h(t, x, y)$, defined as follows. For $x, t\geq0$ and $y\in
\R,$
%
\begin{eqnarray}\label{def:h}
h(t,x,y)&=&\bigl((t+y)\wedge0+x\bigr)_+-(y\wedge 0+x)_+\nonumber
\\&=&\int_{-t}^0\mathbf{1}_{[0, x]}(u-y)\,\mathrm{d}u=
\cases{
t, &\quad$-y<x, y<-t$,\cr
x+t+y,&\quad$-t-y<x\leq-y, y<-t$,\cr
-y, &\quad$-y<x, -t\leq y<0,$\cr
x, &\quad$0<x\leq-y, -t\leq y<0,$\cr
0, &\quad otherwise.
}
\end{eqnarray}
Note that Kaj and Taqqu \cite{KT08}, page 388, interpreted the function
\[
K_t(y,x):=(t-y)_+\wedge x-(-y)_+\wedge x=h(t,x, -x-y)
\]
in the context of the infinite source Poisson model as a function of
the starting time $y$ and the duration $x$ of a session that measures
the length of the time interval contained in $[0, t]$ during which the
session is active. Therefore, $h(t,x,y)=K_t(-x-y,x)$ measures the
length of the time interval contained in $[0, t]$ during which the
session with duration $x$ and finishing time $-y$ is active. Obviously,
by this definition, $h(t,x,y)\not=0$ if and only if $y<0$ and $-y-x<t$,
that is, $x+y>-t$.

In this paper, we define
%
\begin{eqnarray}\label{def:g}
g_s(t,x,y):\!&=&h(t+s,x,y)-h(s,x,y)\nonumber
\\
&=&\int_{-t-s}^0\mathbf{1}_{[0,
x]}(u-y)\,\mathrm{d}u-\int_{-s}^0\mathbf{1}_{[0, x]}(u-y)\,\mathrm{d}u
\\
&=&\int_{-t-s}^{-s}\mathbf{1}_{[0, x]}(u-y)\,\mathrm{d}u=\int_{-t}^0\mathbf{1}_{[0,
x]}(u-y-s)\,\mathrm{d}u=h(t, x, y+s)\nonumber
\end{eqnarray}
for all $t\geq0$, $x\geq0$, $y\in\R$ and any given $s>0$. Then,
according to the definition of an integral of random measure (see
\cite{ST94}, Chapter 3), we can directly verify that for
$H\in(\frac{1}{2}, 1)$,
%
\begin{eqnarray}\label{intr-3}
&&\sqrt{C_H}\int_{0}^{\infty}\!\!\! \int_{\R_-} g_s(t,x, y) x^{H-2}W(\mathrm{d}x,\mathrm{d}
y)\nonumber
\\
&&\quad=\sqrt{C_H}\int_{0}^{\infty}\!\!\! \int_{\R_-} h(t, x, y+s) x^{H-2}W(\mathrm{d} x,\mathrm{d}
y)
\\&&\quad\stackrel{d}{=}\sqrt{C_H}\int_{0}^{\infty}\!\!\! \int_{\R_-} h(t, x, y) x^{H-2}W(\mathrm{d}
x,\mathrm{d}y)\stackrel{d}{=}B^{H}(t),\nonumber
\end{eqnarray}
where $W(\mathrm{d}x,\mathrm{d}y)$ is a Gaussian random measure on
$\mathbb{R}_+\times\mathbb{R}$ with control measure
$\mathrm{d}x\,\mathrm{d}y$ (see Section \ref{2} for its definition),
$C_H=H(2H-1)(1-H)(3-2H)$ and the notation $\stackrel{d}{=}$ denotes
identification of finite-dimensional distributions. In fact, the last
equality is taken from Gaigalas \cite{G06}, page 454, although the
constant $C_H$ is omitted in Gaigalas' representation. Therefore,
%
\begin{equation}\label{s2-3}
C_H\int_0^\infty\!\!\!\int_{\mathbb{R}_-}\frac{g^2_s(t, x,
y)}{x^{4-2H}}\,\mathrm{d}x\,\mathrm{d}y=\frac{1}{2}\bfE[(B^H(t))^2]=\frac{t^{2H}}{2}.
\end{equation}

From the representations (\ref{intr-1}) and (\ref{intr-3}), inspired by
(\ref{intr-11}), it seems reasonable that the law of the process
$B^{H}$ can be approximated by that of some process similar to
$\hat{Y}_{\ep}$ with kernel $g_s(t, x, y)$.

In this paper we define a sequence of processes $\{Y_n(t), t\in[0,
1]\}_{n\geq1}$ as follows:
%
\begin{equation}\label{intr-5}
Y_{n}(t)=n\sqrt{2C_H}\int_0^n\!\!\!\int_{-n}^0 g_s(t,x,y)x^{H-2}
\sqrt{x|y|}(-1)^{N_{n}(x,y)}\,\mathrm{d}x\,\mathrm{d}y
\end{equation}
for $H\in(1/2, 1)$. Here, $\{N_n(x,y), (x,y)\in\R_+ \times\R _-\}$ is a
Poisson process with intensity $n$ (see Definition \ref{d21}).

The main purpose of this paper is to show that the law of $Y_n$
converges to the law of $B^{H}$ for $H\in(1/2, 1)$. We note that the
kernel $g_s(t, x, y)$ in (\ref{intr-5}) cannot be
separated by the arguments $(x, y)$, unlike the kernel function in
(\ref{intr-11}). This difference is not trivial. As we will see in
Remark~\ref{rem3.1}, it causes many real technical difficulties.

The rest of the paper is organized as follows. Section \ref{2} is devoted to
introducing the necessary definitions, notation and the main result. In
Section \ref{3} we prove the family of processes $\{Y_{n}(t)\}_{n\geq1}$
given by (\ref{intr-5}) is tight in $\mathcal{C}[0,1]$. In Section \ref{4} we
prove that as $n$ tends to infinity, the finite-dimensional
distributions of $\{Y_{n}(t)\}$ converge weakly to those of the
fractional Brownian motion $B^{H}$ with $H\in(1/2,1)$ and hence
$\{Y_{n}(t)\}$ converges weakly in $\mathcal{C}[0,1]$ to the fractional
Brownian motion $B^{H}$.

\section{Preliminaries}\label{2}
We now give the definitions of the Brownian sheet and Poisson processes
on $\R\times\R$. Let $\mathscr{F}$ be the Borel algebra on
$\R\times\R$. $\nu$ and $\mu$ denote a $\sigma$-finite measure
and the
Lebesgue measure on $\R\times\R$, respectively.

\begin{defn}\label{d21}
Given a positive constant $\beta>0$, a random set
function $N(\cdot)$ on the measure space
$(\R\times\R, \mathscr{F},\nu)$ is called the \textit{Poisson random
measure with density measure} $\beta\nu$ if it satisfies the following
conditions:
\begin{enumerate}[(iii)]
\item[(i)] for every $A\in\F$ with $\mu(A)<\infty$, $N(A)$ is a Poisson
random variable with parameter $\beta\nu(A)$ defined on the same
probability space;
\item[(ii)] if $A_{1},\ldots, A_n\in\F$ are disjoint
and all have finite measure, then the random variables $N(A_1),\ldots,
N(A_n)$ are independent;
\item[(iii)] if $A_1,A_2,\ldots\in\mathscr{F}$ are disjoint and
$\nu(\bigcup_{i=1}^{\infty} A_{i})<\infty$, then $
N(\bigcup_{i=1}^{\infty}A_i)=\sum _{i=1}^{\infty} N(A_i)$ a.s.
\end{enumerate}
If $N(\cdot)$ is a Poisson random measure with density measure
$\beta\mu$, then we define
\[
N(s,t)=
\cases{
N([0,s]\times[0,t]), &\quad$s\ge0,$ $t\ge0,$\cr
N([0,s]\times[t,0]),&\quad$s\ge0,$ $t\le0,$\cr
N([s,0]\times[0,t]),&\quad$s\le0,$ $t\ge0,$\cr
N([s,0]\times[t,0]), &\quad$s\le0,$ $t\le0$
}
\]
and call $N=\{N(s, t), (s, t)\in\R\times\R\}$ the two-parameter
Poisson process with intensity $\beta$ in $\R\times\R$.
\end{defn}

It is not hard to see that $N=\{N(s, t), (s,
t)\in\mathbb{R}_+\times\mathbb{R}_-\}$ is independent with $N_1=\{N(s,
t), (s, t) \in\mathbb{R}_+\times\mathbb{R}_+\}$, which is the ordinary
two-parameter process in $\mathbb{R}_+\times\mathbb{R}_+$, and for any
$(s,t)\in\R_+\times\R_-$, $\{N(s,t)\}$ has the same finite-dimensional
distributions as those of $\{N(s, |t|)\}.$

\begin{defn}
Consider a random set function $W(\cdot)$ on the measure space
$(\R\times\R, \F,\nu)$ such that:
\begin{enumerate}[(iii)]
\item[(i)] for every $A\in\F$ with $\nu(A)<\infty$, $W(A)$ is a
centered Gaussian random variable defined on the same probability space
with variance $2\nu(A)$;
\item[(ii)] if $A_{1},\ldots, A_{n}\in\F$ are
disjoint and have finite measure, then the random variables
$W(A_{1}),\ldots,W(A_{n})$ are independent;
\item[(iii)] if $A_1,A_2, \ldots\in\mathscr{F}$ are disjoint and
$\nu(\bigcup_{i=1}^{\infty}A_i)<\infty$, then $W(\bigcup_{i=1}^{\infty}
A_i)=\sum_{i=1}^{\infty} W(A_i)$ a.s.
\end{enumerate}
We then call $W(\cdot)$ a \textit{Gaussian random measure on
$\R\times\R$ with control measure} $\nu$. In particular, if
$W(\cdot)$
is a Gaussian random measure on $\R\times\R$ with control measure
$\frac{1}{2}\mu$, then we define
\[
B(s,t)=
\cases{
W([0,s]\times[0,t]), &\quad$s\geq0,$ $t\geq0,$\cr
W([0,s]\times[t,0]),&\quad$s\geq0,$ $t\leq0,$\cr
W([s,0]\times[0,t]),&\quad$s\leq0,$ $t\geq0,$\cr
W([s,0]\times[t,0]),&\quad$s\leq0,$ $t\leq0$
}
\]
and call $B=\{B(s, t), (s, t)\in\R\times\R\}$ the \textit{two-parameter
Brownian sheet in} $\R\times\R$.
\end{defn}

Similarly, we have that $B=\{B(s, t), (s,
t)\in\mathbb{R}_+\times\mathbb{R}_-\}$ is independent of $B_1=\{B(s,
t), (s, t) \in\mathbb{R}_+\times\mathbb{R}_+\}$, which is the ordinary
Brownian sheet in $\mathbb{R}_+\times\mathbb{R}_+$, and for any
$(s,t)\in\R\times\R$, we have $B(s,t)\stackrel{d}{=}B(|s|,|t|).$ Hence,
from (\ref{intr-3}) it is easy to check that
%
\begin{eqnarray}\label{pre-1}
B^{H}(t) \stackrel{d}{=}\sqrt{2C_H} \int_{0}^{\infty}\!\!\!\int
_{\mathbb{R}_-} g_s(t,x, y) x^{H-2} B(\mathrm{d} x,\mathrm{d} y).
\end{eqnarray}

Let $\operatorname{sgn}(x)=1$ if $x\geq0$ and $\operatorname{sgn}(x)=-1$ if $x< 0$. We
have the following conclusion, which essentially parallels
Bardina and Jolis \cite{BJ001}, Theorem 1.1. The proof is omitted.

\begin{lem}\label{pre-lem-1}
Suppose that $N=\{N(s,t), (s, t)\in\R\times\R\}$ is a two-parameter
Poisson process with intensity $1$ in $\R\times\R$. For any $S>0$ and
$T>0$, let
%
\begin{eqnarray}\label{pre-2}
B_n(u,v)=\operatorname{sgn}(uv) n\int_0^u\!\!\!\int_0^v\sqrt{|xy|}(-1)^{N(x\sqrt
{n},y\sqrt{n})}\,\mathrm{d}x\,\mathrm{d}y
\end{eqnarray}
for any $|u|\leq S$, $|v|\leq T$. The finite-dimensional distributions
of $B_n$ then converge weakly to those of a two-parameter Brownian
sheet $B=\{B(u,v), |u|\leq S, |v|\leq T\}$.
\end{lem}

Naturally, (\ref{pre-1}) and (\ref{pre-2}) suggest that we consider the
following approximation of $B^H$ for $H\in(1/2, 1)$:
%
\begin{eqnarray}\label{pre-3}
Y_{n}(t)=n\sqrt{2C_H}
\int_0^n\!\!\!\int_{-n}^0g_s(t,x,y)x^{H-2}\sqrt{x|y|}(-1)^{N_n(x,y)}\,\mathrm{d}x\,\mathrm{d}y
\end{eqnarray}
for $n\in\mathbb{N}=\{1,2,\ldots\}$ and $t\in[0,1]$, where
$N_n=\{N_n(x, y)\}=\{N(x\sqrt{n},y\sqrt{n})\}$ is the two-parameter
Poisson process with intensity $n$.

The main result of this paper is as follows.

\begin{thhm}\label{s4-thhm1}
For any fixed $s>0$, the law of the process $\{Y_{n}(t),t\in[0,1]\}$
given by (\ref{pre-3}) converges weakly to the law of $\{B^{H}(t),
t\in[0,1]\}$ for $H\in(\frac{1}{2},1)$ in $\mathcal{C}[0,1]$.
\end{thhm}

\begin{rem}\label{r21}
If we define $\bar{h}(t, x, y):=h(t,x,-y)$, then $\bar{h}(t,x,y)$ has a
more natural physical interpretation. It measures the length of the
time interval contained in $[0, t]$ during which the session with
duration $x$ and finishing time $y$ is active. Define
\[
\bar{Y}_{n}(t):=n\sqrt{2C_H}
\int_0^n\!\!\!\int_0^n\bar{g}_s(t,x,y)x^{H-2}\sqrt{xy}(-1)^{N_n(x,
y)}\,\mathrm{d}x\,\mathrm{d}y,
\]
where $\bar{g}_s(t,x,y)=\bar{h}(t+s,x,y)-\bar{h}(s, x, y)=\bar {h}(t,
x, y-s)$. It is then easy to see that Theorem~\ref{s4-thhm1} holds for
$\bar{Y}_n$.
\end{rem}

Note that $s$ is a given positive number. In the sequel we will treat
it as a constant. In addition, since the parameter $C_H$ cannot affect
our discussion, we will take it to be $1$ in order to simplify
matters.

We now introduce some auxiliary notation. In the sequel, $\infty$ and
$-\infty$ will denote positive infinity and negative infinity,
respectively. For all $x'\geq x\geq0$, $y'\geq y\geq0$, by Cairoli and
Walsh \cite{CW75}, we define
%
\begin{equation}\label{intr-6}
\Delta_{(x,y)}N(x',y')=N\bigl((x,x']\times(y,y']\bigr).
\end{equation}
For any $t\in[0, 1]$, $x>0$ and $y\leq0$, we let
\[
\phi_t(x, y):=g_s(t,x,y)/x^{2-H}\geq0.
\]
Define a function $F_{n, f}(x, y)$ as follows:
%
\begin{equation}\label{intr-7}
F_{n, f}(x, y)=\sqrt{x|y|}n^2\int_x^{\infty}\!\!\!\int_y^0 f(x_2,y_2)
\sqrt{x_2|y_2|}\e^{-2n[x(y_2-y)-(x_2-x)y_2]}\,\mathrm{d}x_2\,\mathrm{d}y_2,
\end{equation}
where $n>0$, $x\geq0, y\leq0$ and $f$ is a measurable function such
that the integral is meaningful. Obviously, if $0\leq f\leq g$, then
$0\leq F_{n, f}(x, y)\leq F_{n, g}(x, y).$

\section{Tightness of $Y_n$ in ${\mathcal C}[0, 1]$}\label{3}

Let $Y_n=\{Y_{n}(t),t\in[0,1]\}$ be the process defined by
(\ref{pre-3}). The purpose of this section is to prove the tightness of
the processes $\{Y_n\}_{n\geq1}$.

\begin{prop}\label{s3-prop1}
The family of laws of the processes $\{Y_n\}_{n\geq1}$ is tight in
${\mathcal C}[0, 1]$.
\end{prop}

To prove the proposition, we need the following lemmas.

\begin{lem}\label{s3-lem1}
Let $\Omega_1=\{0<x_1\leq x_2, y_2\leq y_1\leq0\}$ and
$\Omega_2=\{0<x_1\leq x_2, y_1<y_2\leq0\}$. For a non-negative function
$f(x, y)$, if $\int_0^\infty\!\!\int_{-\infty}^0 f^2(x,
y)\,\mathrm{d}x\,\mathrm{d}y<\infty$, then for any $u, v>0$,
%
\begin{eqnarray}\label{s3-lem1-0}
&&\bfE\biggl[\biggl(n\int_0^u\!\!\!\int_{-v}^0 f(x,y)\sqrt
{x|y|}(-1)^{N_{n}(x,y)}\,\mathrm{d}x\,\mathrm{d}y\biggr)^{2}\biggr] \nonumber
\\[-8pt]\\[-8pt]
&&\quad \leq2\bigl(I_1(n,
f)+I_2(n, f)\bigr), \nonumber
\end{eqnarray}
where
%
\begin{eqnarray}
\label{s3-1}I_1(n, f)
&=&
n^{2}\int_{\Omega_1}\prod_{i=1}^{2}\bigl(f(x_i,y_i)\sqrt{x_i|y_i|}\bigr)\e^{-2n(x_1y_1-x_2y_2)}\,\mathrm{d}x_1\,\mathrm{d}y_1\,\mathrm{d}x_2\,\mathrm{d}y_2,
\\
\label{s3-2}I_2(n, f)
&=&
n^{2}\int_{\Omega_2}\prod_{i=1}^{2}\bigl(f(x_i,y_i)\sqrt{x_i|y_i|}\bigr)\e^{-2n[x_1(y_2-y_1)-(x_2-x_1)y_2]}\,\mathrm{d}x_1\,\mathrm{d}y_1\,\mathrm{d}x_2\,\mathrm{d}y_2.\quad
\end{eqnarray}
\end{lem}

\begin{pf}
Let $I(n,x_{1,2},y_{1,2})=\bfE[
(-1)^{N_{n}(x_{1},y_{1})+N(x_{2},y_{2})}]$. By Fubini's theorem, the
left-hand side of (\ref{s3-lem1-0}) is equal to
%
\begin{equation}\label{s3-l1-1}
n^{2}\int_0^u\!\!\!\int_{-v}^0\int_0^u\!\!\!\int
_{-v}^0f(x_1,y_1)f(x_2,y_2)\sqrt{x_1|y_1|}\sqrt{x_2|y_2|}
I(n,x_{1,2},y_{1,2})\,\mathrm{d}x_1\,\mathrm{d}y_1\,\mathrm{d}x_2\,\mathrm{d}y_2.
\end{equation}
Define $\Omega_3=\{x_1>x_2>0, 0\geq y_1\geq y_2\}$ and $\Omega_4=\{x_1>
x_2>0, y_1<y_2\leq0\}$. Then, by (\ref{intr-6}), $I(n,
x_{1,2},y_{1,2})$ equals
%
\begin{eqnarray}
&&\bfE\bigl[
(-1)^{N_{n}(x_{1},|y_{1}|)+N_{n}(x_{2},|y_{2}|)}\bigr]\nonumber
\\[-8pt]\\[-8pt]
&&\quad=\bfE\bigl[(-1)^{\Delta_{(0,0)}N_{n}(x_{1},|y_{1}|)
+\Delta_{(0,0)}N_{n}(x_{2},|y_{2}|)}\bigr].\nonumber
\end{eqnarray}
Note that $\sum_{i=1}^{2}\Delta_{(0,0)}N_{n}(x_{i},|y_i|)$ is\vadjust{\goodbreak} equal to
the sum of the increments of the Poisson process over some disjoint
rectangles, and the rectangles which contribute to the value of
$I(n,x_{1,2},y_{1,2})$ are those which appear only once. Since
two-parameter Poisson processes have independent increments, after some
simple calculation, we obtain that on $\Omega_1$,
\begin{eqnarray}\label{s3-l1-2}
&&
I(n,x_{1,2},y_{1,2})\nonumber
\\
&&\quad=
\bfE\bigl[(-1)^{\Delta_{(0,|y_{1}|)}N_{n}(x_{1},|y_{2}|)+\Delta_{(x_{1},|y_{1}|)}N_{n}(x_{2},|y_{2}|)+\Delta_{(x_{1},0)}N_{n}(x_{2},|y_{1}|)}\bigr]\nonumber
\\[-8pt]\\[-8pt]
&&\quad=
\bfE\bigl[(-1)^{\Delta_{(0,|y_{1}|)}N_{n}(x_{1},|y_{2}|)}\bigr]\bfE\bigl[(-1)^{\Delta_{(x_{1},|y_{1}|)}
N_{n}(x_{2},|y_{2}|)}\bigr]\bfE\bigl[(-1)^{\Delta_{(x_{1},0)}N_{n}(x_{2},|y_{1}|)} \bigr]\nonumber
\\
&&\quad=
\exp\{-2n(x_1y_1-x_2y_2)\}.\nonumber
\end{eqnarray}
Using the same method as above, we obtain that
%
\begin{eqnarray}
\label{s3-l1-3}
I(n, x_{1,2}, y_{1,2})&=&\exp\{-2n
[x_1(y_2-y_1)-(x_2-x_1)y_2]\}\qquad\mbox{on }\Omega_2,
\\
\label{s3-l1-4}
I(n, x_{1,2},y_{1,2})&=&\exp\{-2n
[x_2(y_1-y_2)-(x_1-x_2)y_1]\}\qquad\mbox{on }\Omega_3,
\\
\label{s3-l1-5}
I(n, x_{1,2}, y_{1,2})&=&\exp\{-2n(x_2y_2-x_1y_1)\}\qquad\mbox{on }\Omega_4.
\end{eqnarray}
Substituting (\ref{s3-l1-2})--(\ref{s3-l1-5}) into (\ref{s3-l1-1}) and
using a change of the integration variables if necessary, we can easily
obtain (\ref{s3-lem1-0}).
\end{pf}

\begin{lem}\label{s3-lem2} If $\int_0^\infty\!\!\int_{-\infty}^0
f^2(x, y)\,\mathrm{d}x\,\mathrm{d}y<\infty$, then $I_1(n, f)$ defined by
(\ref{s3-1}) is such that
%
\begin{eqnarray}\label{s3-l2-0}
I_1(n,
f)\leq\frac{37}{8}\int_0^\infty\!\!\!\int_{-\infty}^0f^2(x,y)\,\mathrm{d}x\,\mathrm{d}y.
\end{eqnarray}
\end{lem}

\begin{pf}
Define $B=\{0<x_1\leq x_2\leq2x_1\}$, $ C=\{2y_1\leq
y_2\leq y_1\leq0\}$ and $A=B\cap C$. Let $I_1^{A}(n)$, $I_1^{B}(n)$ and
$I_1^{C}(n)$ denote the integral (\ref{s3-1}), where $\Omega_1$ is
replaced by $A$, $\Omega_1\setminus B$ and $\Omega_1\setminus C$,
respectively. Then,
%
\begin{eqnarray}\label{s3-l2-1}
I_1(n, f)\leq I_1^{A}(n, f)+I_1^{B}(n, f)+I_1^{C}(n, f).
\end{eqnarray}

Using the elementary inequality $2ab\leq a^{2}+b^{2},$ from
(\ref{s3-1}), we have that
%
\begin{eqnarray}\label{s3-lem1-21}
I_1^{A}(n, f)\leq\tfrac{1}{2}\bigl(I_{11}(n)+I_{12}(n)\bigr),
\end{eqnarray}
where $I_{11}(n)$ is
\begin{eqnarray*}
&&
n^{2}\int_Af^2(x_1,y_1)x_1|y_1|\e^{-2n(x_1y_1-x_2y_2)}\,\mathrm{d}x_1\,\mathrm{d}y_1\,\mathrm{d}x_2\,\mathrm{d}y_2
\\
&&\quad=
n^{2}\int_{0}^{\infty}\!\!\!\int_{-\infty}^{0}\int_{x_1}^{2x_1}\!\!\!\int_{2y_1}^{y_1}
f^2(x_1,y_1)x_1|y_1|\e^{-2n(x_1y_1-x_2y_2)}\,\mathrm{d}x_1\,\mathrm{d}y_1\,\mathrm{d}x_2\,\mathrm{d}y_2
\end{eqnarray*}
and $I_{12}(n)$ is
\begin{eqnarray*}
&&
n^{2}\int_Af^2(x_2,y_2)
x_2|y_2|\e^{-2n(x_1y_1-x_2y_2)}\,\mathrm{d}x_1\,\mathrm{d}y_1\,\mathrm{d}x_2\,\mathrm{d}y_2
\\
&&\quad=
n^{2}\int_{0}^{\infty}\!\!\!\int_{-\infty}^{0}\int_{x_1}^{2x_1}\!\!\!\int_{2y_1}^{y_1}f^2(x_2,y_2)
x_2|y_2|\e^{-2n(x_1y_1-x_2y_2)}\,\mathrm{d}x_1\,\mathrm{d}y_1\,\mathrm{d}x_2\,\mathrm{d}y_2.
\end{eqnarray*}
Since, for any $(x_1,y_1,x_2,y_2)\in\Omega_1$,
%
\begin{equation}\label{s3-lem1-a0}
x_1y_1-x_2y_2\geq(x_{2}-x_{1})|y_1|+(|y_2|-|y_1|)x_1,
\end{equation}
$I_{11}(n)$ is bounded by
\[
n^{2}\int_{0}^{\infty}\!\!\!\int_{-\infty}^{0}\int_{x_1}^{2x_1}\!\!\!\int
_{2y_1}^{y_1}f^2(x_1,y_1) x_{1}|y_1|
\e^{-2n[(x_{2}-x_{1})|y_1|+(|y_2|-|y_1|)x_1]}\,\mathrm{d}x_1\,\mathrm{d}y_1\,\mathrm{d}x_2\,\mathrm{d}y_2.
\]
Integrating with respect to $x_{2}$ and then with respect to $y_{2}$ in
the above integral, we obtain the following bound:
%
\begin{equation}\label{s3-lem1-a2}
I_{11}(n)\leq\frac{1}{4} \int_0^\infty\!\!\!\int_{-\infty}^0
f^2(x_1,y_1)\,\mathrm{d}x_1\,\mathrm{d}y_1.
\end{equation}
Furthermore, in the region $A$, it follows from (\ref{s3-lem1-a0}) that
\[
x_1y_1-x_2y_2\geq(x_{2}-x_{1})|y_2|/2+(|y_2|-|y_1|)x_2/2.
\]
By the same argument as above, we then have that
%
\begin{equation}\label{s3-lem1-a3}
I_{12}(n)\leq\int_0^\infty\!\!\!\int_{-\infty}^0
f^2(x_2,y_2)\,\mathrm{d}x_2\,\mathrm{d}y_2.
\end{equation}
Therefore, from (\ref{s3-lem1-21}), (\ref{s3-lem1-a2}) and
(\ref{s3-lem1-a3}), it follows that
%
\begin{equation}\label{s3-l2-2}
I_1^{A}(n, f)\leq\frac{5}{8}\int_0^\infty\!\!\!\int_{-\infty}^0
f^2(x,y)\,\mathrm{d}x\,\mathrm{d}y.
\end{equation}

We now consider $I_1^{B}(n, f)$. Note that
\begin{eqnarray*}
 2(x_1y_1-x_2y_2)
 &=&
 2(x_2-x_1)|y_1|+2(|y_2|-|y_1|)x_2
\\
&\geq&
(x_{2}-x_{1})|y_{1}|+(|y_{2}|-|y_{1}|)x_{1}+ x_{1}|y_{1}|
\\
&=&
(x_2-x_1)|y_1|+|y_2|x_1
\\
&\geq&
\tfrac{1}{2}x_2|y_1|+|y_2|x_1
\end{eqnarray*}\vadjust{\goodbreak}%
for any $(x_1,y_1,x_2,y_2)\in\Omega_1\setminus B=\{x_2> 2x_1>0,
y_2\leq
y_1\leq0\}$. Define $\tilde{I}_{11}(n)$ and $\tilde{I}_{12}(n)$ as
$I_{11}(n)$ and $I_{12}(n)$, respectively, with $\Omega_1\setminus B$
instead of $A$. Then,
\begin{eqnarray*}
\tilde{I}_{11}(n)
&\leq&
n^{2}\int_0^\infty\!\!\!\int_{-\infty}^0\int_{2x_1}^\infty\!\!\int_{-\infty}^{y_1}
f^2(x_1,y_1)x_1|y_1|\e^{-n(|y_2|x_1+(1/2)x_2|y_1|)}\,\mathrm{d}x_1\,\mathrm{d}y_1\,\mathrm{d}x_2\,\mathrm{d}y_2,
\\
\tilde{I}_{12}(n)
&\leq&
n^{2}\int_0^\infty\!\!\!\int_{-\infty}^0\int_{2x_1}^\infty\!\!\int_{-\infty}^{y_1}
f^2(x_2,y_2)x_{2}|y_{2}|\e^{-n(|y_2|x_1+(1/2)x_2|y_1|)}\,\mathrm{d}x_1\,\mathrm{d}y_1\,\mathrm{d}x_2\,\mathrm{d}y_2.
\end{eqnarray*}
Integrating in $\tilde{I}_{11}(n)$ with respect to $x_2$ and $y_2$, and
in $\tilde{I}_{12}(n)$ with respect to $x_1$ and $y_1$, respectively,
we obtain that
\[
\tilde{I}_{11}(n)+\tilde{I}_{12}(n)\leq4\int_0^\infty\!\!\!\int
_{-\infty}^0f^2(x,y)\,\mathrm{d}x\,\mathrm{d}y.
\]
Since $I_1^{B}(n, f)\leq\frac{1}{2}(
\tilde{I}_{11}(n)+\tilde{I}_{12}(n))$, we have
%
\begin{equation}\label{s3-l2-3}
I_1^{B}(n, f)\leq2\int_0^\infty\!\!\!\int_{-\infty}^0
f^2(x,y)\,\mathrm{d}x\,\mathrm{d}y.
\end{equation}

Similarly, for any $(x_1,y_1,x_2,y_2)\in\Omega_1\setminus C=\{x_2\geq
x_1>0, y_2<2y_1\leq0\}$, we have that
\begin{eqnarray*}
2(x_1y_1-x_2y_2)
&=&
2(x_2|y_2|-|y_1|x_1)\geq2(|y_2|-|y_1|)x_2
\\
&\geq&
(|y_2|-|y_1|)x_2+|y_1|x_2
\geq(|y_2|-|y_1|)x_1+|y_1|x_2
\\
&\geq&
\tfrac{1}{2}|y_2|x_1+x_2|y_1|.
\end{eqnarray*}
Therefore, using a similar approach to the one above, we have
%
\begin{eqnarray}\label{s3-l2-4}
I_1^{C}(n, f)\leq2\int_0^\infty\!\!\!\int_{-\infty}^0
f^2(x,y)\,\mathrm{d}x\,\mathrm{d}y.
\end{eqnarray}
Combining (\ref{s3-l2-1}) with (\ref{s3-l2-2})--(\ref{s3-l2-4}), we
get (\ref{s3-l2-0}).
\end{pf}

From (\ref{s2-3}) and Lemma \ref{s3-lem2}, we can immediately get the
following corollary.

\begin{cor}\label{s3-c1}
For each $n>0$, $I_1(n, \phi_t)\leq\frac
{37}{16}t^{2H}$.
\end{cor}

Observe that for every non-negative measurable function $f(x, y)$, by
(\ref{intr-7}) and
(\ref{s3-2}),
%
\begin{equation}
I_2(n, f)=\int_0^\infty\!\! \mathrm{d}x_1\int_{-\infty}^0 f(x_1,y_1)F_{n,
f}(x_1, y_1)\,\mathrm{d}y_1.
\end{equation}
The following lemmas concern $I_2(n, f)$. We focus on the case where
$f=\phi_t$.

\begin{lem}\label{s3-lema}
For any $0\leq t\leq1$, let $S=\{x>0, -t-s<y<-s\}$ and $\tilde{\phi}_t(x,
y)=\phi_{t}(x, y)\mathbf{1}_{S}(x, y) $.
There then exist a non-negative function $\tilde{\Psi}_t(x, y)$ on $\{
x>0, y\leq0\}$ such that for all $n>0$,
%
\begin{equation}\label{s3-la-01}
F_{n,\tilde{\phi}_t}(x_1, y_1)\mathbf{1}_{S}(x_1, y_1)\leq\tilde{\Psi
}_t(x_1, y_1),
\end{equation}
and a positive constant $K_1$ which depends only on $H$, such
that for all $n>0$,
%
\begin{equation}\label{s3-la-02}
I_2(n,\tilde{\phi}_t)\leq\int_{S}\phi_t(x_1, y_1)\tilde{\Psi}_t(x_1,
y_1)\,\mathrm{d}x_1\,\mathrm{d}y_1\leq K_1\sqrt{\frac{s+1}{s}}t^{2H}.
\end{equation}
\end{lem}

\begin{pf}
By (\ref{def:h}) and (\ref{def:g}),
%
\begin{equation}\label{s3-la-1}
\phi_t(x, y)=
\cases{
 t/x^{2-H}, &\quad$-y-s<x,$ $y<-t-s,$\cr
(x+t+s+y)/x^{2-H},&\quad$-t-s-y<x\leq-y-s,$ $y<-t-s,$\cr
-(y+s)/x^{2-H}, &\quad$-y-s<x,$ $-t-s\leq y<-s,$\cr
1/x^{1-H}, &\quad$ 0<x\leq-y-s,$ $-t-s\leq y<-s,$\cr
0, &\quad otherwise.
}
\end{equation}
Let $S_1=S\cap\{0<x<-y-s\}$ and $S_2=S\cap\{x>-y-s\}$. Then,
%
\begin{equation}\label{s3-la-2}
F_{n,\tilde{\phi}_t}(x_1, y_1)\mathbf{1}_{S_1}(x_1, y_1)=\mathbf{
1}_{S_1}(x_1, y_1)\sqrt{x_1}(I_{21}+I_{22}+I_{23}),
\end{equation}
where
\begin{eqnarray*}
I_{21}
&=&
n^2\int_{y_1}^{-s-x_1}\!\!\int_{x_1}^{-s-y_2}\phi
_t(x_2,y_2)\sqrt{x_2|y_1y_2|}
\e^{-2n[x_1(y_2-y_1)-(x_2-x_1)y_2]}\,\mathrm{d}x_2\,\mathrm{d}y_2,
\\
I_{22}
&=&
n^2\int_{y_1}^{-s-x_1}\!\!\int_{-s-y_2}^{\infty}\phi
_t(x_2,y_2)\sqrt{x_2|y_1y_2|}
\e^{-2n[x_1(y_2-y_1)-(x_2-x_1)y_2]}\,\mathrm{d}x_2\,\mathrm{d}y_2,
\\
I_{23}
&=&
n^2\int_{-s-x_1}^{-s}\!\int_{x_1}^{\infty}\phi
_t(x_2,y_2)\sqrt{x_2|y_1y_2|}
\e^{-2n[x_1(y_2-y_1)-(x_2-x_1)y_2]}\,\mathrm{d}x_2\,\mathrm{d}y_2.
\end{eqnarray*}
For $(x_1, y_1)\in S_1$, from (\ref{s3-la-1}), we have that
%
\begin{eqnarray}\label{s3-la-3}
I_{21}
&=&
n^2\int_{y_1}^{-s-x_1}\!\int_{x_1}^{-s-y_2}\frac
{1}{x_2^{1/2-H}}\sqrt{|y_1y_2|}
\e^{-2n[x_1(y_2-y_1)-(x_2-x_1)y_2]}\,\mathrm{d}x_2\,\mathrm{d}y_2
\nonumber
\\
&\leq&
\sqrt{\frac{t+s}{s}}n^2\int_{y_1}^{-s-x_1}\!\int
_{x_1}^{-s-y_2}\frac{1}{x_2^{1/2-H}}|y_2|
\e^{-2n[x_1(y_2-y_1)-(x_2-x_1)y_2]}\,\mathrm{d}x_2\,\mathrm{d}y_2\nonumber
\\[-8pt]\\[-8pt]
&\leq&
t^{H-1/2}\sqrt{\frac{t+s}{s}}\frac{n}{2}\int
_{y_1}^{-s-x_1}\e^{-2nx_1(y_2-y_1)}\,\mathrm{d}\nonumber
y_2
\\
&\leq&
\frac{t^{H-1/2}}{4x_1}C_s=:\psi_1(x_1,y_1),\nonumber
\end{eqnarray}
where $C_s=\sqrt{(1+s)/s}$ is a constant. Similarly, for $(x_1,
y_1)\in
S_1$,
\begin{eqnarray}
I_{22}
&=&
n^2\int_{y_1}^{-s-x_1}\!\int_{-s-y_2}^{\infty}\frac
{-s-y_2}{x_2^{3/2-H}}\sqrt{|y_1y_2|}
\e^{-2n[x_1(y_2-y_1)-(x_2-x_1)y_2]}\,\mathrm{d}x_2\,\mathrm{d}y_2
\nonumber
\\
&\leq&
C_sn^2\int_{y_1}^{-s-x_1}\!\int_{-s-y_2}^{\infty
}(-s-y_2)^{H-1/2}|y_2|
\e^{-2n[x_1(y_2-y_1)-(x_2-x_1)y_2]}\,\mathrm{d}x_2\,\mathrm{d}y_2\nonumber
\\[-8pt]\\[-8pt]
&\leq&
(-s-y_1)^{H-1/2}C_s\frac{n}{2}\int_{y_1}^{-s-x_1}
\e^{-2nx_1(y_2-y_1)}\,\mathrm{d} y_2
\nonumber
\\
&\leq&
\frac{(-s-y_1)^{H-1/2}}{4x_1}C_s=:\psi_2(x_1,y_1)\nonumber
\end{eqnarray}
and
%
\begin{eqnarray}\label{s3-la-4}
I_{23}
&=&
n^2\int_{-s-x_1}^{-s}\int_{x_1}^{\infty}\frac
{-s-y_2}{x_2^{3/2-H}}\sqrt{|y_1||y_2|}
\e^{-2n[x_1(y_2-y_1)-(x_2-x_1)y_2]}\,\mathrm{d}x_2\,\mathrm{d}y_2\nonumber
\\
&\leq&
C_sn^2\int_{-s-x_1}^{-s}\int_{x_1}^{\infty}\frac
{x_1}{x_1^{3/2-H}}|y_2|
\e^{-2n[x_1(y_2-y_1)-(x_2-x_1)y_2]}\,\mathrm{d}x_2\,\mathrm{d}y_2
\\
&\leq&
C_s\frac{1}{4x_1^{3/2-H}}=:\psi_3(x_1,y_1).\nonumber 
\end{eqnarray}
In addition, for $(x_1,y_1)\in S_2$,
\begin{eqnarray}\label{s3-la-7}
F_{n,\tilde{\phi}_t}(x_1,y_1)
&=&
n^2\int_{y_1}^{-s}\!\!\!\int_{x_1}^{\infty}\frac
{-s-y_2}{x_2^{3/2-H}}\sqrt{x_1|y_1y_2|}
\e^{-2n[x_1(y_2-y_1)-(x_2-x_1)y_2]}\,\mathrm{d}x_2\,\mathrm{d}y_2
\nonumber
\\[-8pt]\\[-8pt]
&\leq&
C_s\frac{-s-y_1}{4x_1^{2-H}}=:\psi
_4(x_1,y_1).\nonumber
\end{eqnarray}
Define $\tilde{\Psi}_t(x_1, y_1)$ as
\[
\sqrt{x_1}\bigl(\psi_1(x_1, y_1)+\psi_2(x_1, y_1)+\psi_3(x_1,
y_1)\bigr)\mathbf{1}_{S_1}(x_1, y_1)+\psi_4(x_1, y_1)\mathbf{1}_{S_2}(x_1,y_1).
\]
Obviously, $\tilde{\Psi}_t(x_1,y_1)$ is positive and (\ref{s3-la-01})
follows from (\ref{s3-la-2})--(\ref{s3-la-7}). Note that
\begin{eqnarray*}
I_2(n,\tilde{\phi}_t)&=&\int_S\tilde{\phi}_t(x_1,
y_1)F_{n,\tilde{\phi}_t} (x_1,
y_1)\,\mathrm{d}x_1\,\mathrm{d}y_1\\
&\leq&\int_S\tilde{\phi}_t(x_1,
y_1)\tilde{\Psi }_t(x_1,y_1)\,\mathrm{d}x_1\,\mathrm{d}y_1.
\end{eqnarray*}
With some basic calculations, we obtain the following results:
\begin{eqnarray*}
\int_{S_1}\phi_t(x, y)\sqrt{x_1}\psi_1(x_1,y_1)\,\mathrm{d}x_1\,\mathrm{d}y_1
&=&
\frac{C_s}{4}\int_{-t-s}^{-s}\int_0^{-y_1-s}\frac{t^{H-1/2}}{x_1^{3/2-H}}\,\mathrm{d}x_1\,\mathrm{d}y_1=\frac{C_st^{2H}}{4H^2-1};
\\
\int_{S_1}\phi_t(x, y)\sqrt{x_1}\psi_2(x_1,y_1)\,\mathrm{d}x_1\,\mathrm{d}y_1
&=&
\frac{C_s}{4}\int_{-t-s}^{-s}\int_0^{-y_1-s}\frac{(-s-y_1)^{H-1/2}}{x_1^{3/2-H}}\,\mathrm{d}x_1\,\mathrm{d}y_1
\\
&=&
\frac{C_st^{2H}}{4(2H-1)H};
\\
\int_{S_1}\phi_t(x, y)\sqrt{x_1}\psi_3(x_1,y_1)\,\mathrm{d}x_1\,\mathrm{d}y_1
&=&
\frac{C_s}{4}\int_{-t-s}^{-s}\int_0^{-y_1-s}\frac{1}{x_1^{2-2H}}\,\mathrm{d}x_1\,\mathrm{d}y_1 =\frac{C_st^{2H}}{8(2H-1)H};
\\
\int_{S_2}\phi_t(x, y)\psi_4(x_1,y_1)\,\mathrm{d}x_1\,\mathrm{d}y_1
&=&
\frac{C_s}{4}\int_{-t-s}^{-s}\int_{-y_1-s}^\infty\frac{(-y_1-s)^2}{x_1^{4-2H}}\,\mathrm{d}x_1\,\mathrm{d}y_1
\\
&=&\frac{C_st^{2H}}{8(3-2H)H}.
\end{eqnarray*}
From the above integrals, (\ref{s3-la-02}) follows with
$K_1=\frac{3}{8H(2H-1)}+\frac{1}{8H(3-2H)}+\frac{1}{4H^2-1}$.
\end{pf}

\begin{lem}\label{s3-lem4}
For $0\leq t\leq1$ and each constant $M<-t-s$, let $G(M)=\{-t-s-y<x,
M<y<-t-s\}$ and $\hat{\phi}_t(x, y)=\phi_t(x,y)\mathbf{1}_{G(M)}(x,y)$.
There then exist a non-negative function $\hat{\Psi}_t(x, y)$ on $\{
x>0,y\leq0\}$
such that for all $n>0$,
%
\begin{equation}\label{s3-l4-01}
F_{n,\hat{\phi}_t}(x_1, y_1)\mathbf{1}_{G(M)}(x_1,y_1)\leq\hat{\Psi
}_t(x_1, y_1),
\end{equation}
and a positive constant $K_2$ which only depends on $H$, such that
%
\begin{equation}\label{s3-l4-02}
I_2(n, \hat{\phi}_t)\leq\int_{G(M)}\phi_t(x_1, y_1)\hat{\Psi}_t(x_1,
y_1)\,\mathrm{d}x_1\,\mathrm{d}y_1\leq K_2\sqrt{\frac{|M|}{s}}t^{2H}.
\end{equation}
\end{lem}

\begin{pf}
From (\ref{s3-la-1}), we have that
%
\begin{equation}\label{s3-l4-2}
F_{n,\hat{\phi}_t}(x_1, y_1)=\hat{I}_{21}+\hat{I}_{22},
\end{equation}
where
\begin{eqnarray*}
\hat{I}_{21}
&=&
n^2\int_{y_1}^{-t-s}\!\!\!\int_{-s-y_2}^{\infty}\frac{t}{x_2^{3/2-H}}\sqrt{x_1|y_1y_2|}\e^{-2n[x_1(y_2-y_1)-(x_2-x_1)y_2]}\,\mathrm{d}x_2\,\mathrm{d}y_2,
\\
\hat{I}_{22}
&=&
n^2\int_{y_1}^{-t-s}\!\!\!\int_{-x_1}^{-s-y_2}\frac{x_2+t+s+y_2}{x_2^{3/2-H}}\sqrt{x_1|y_1y_2|}\e^{-2n[x_1(y_2-y_1)-(x_2-x_1)y_2]}\,\mathrm{d}x_2\,\mathrm{d}y_2.
\end{eqnarray*}
For $(x_1,y_1)\in G(M)$, since $|y_1|<|M|$,
\begin{eqnarray}\label{s3-l4-3}
\hat{I}_{21}\vee\hat{I}_{22}
&\leq&
\sqrt{\frac{|M|}{s}}\frac{t}{x_1^{1-H}}\int_{y_1}^{-t-s}\!\!\mathrm{d}y_1
\int_{-s-y_2}^{\infty}n^2|y_2|\e^{-2n[x_1(y_2-y_1)-(x_2-x_1)y_2]}\,\mathrm{d}x_2\,\mathrm{d}y_2
\nonumber
\\[-8pt]\\[-8pt]
&\leq&
\sqrt{\frac{|M|}{s}}\frac{t}{4x_1^{2-H}}.\nonumber
\end{eqnarray}
Therefore, from (\ref{s3-l4-2}) and (\ref{s3-l4-3}), it follows that
\[
F_{n,\hat{\phi}_t}(x_1, y_1)\mathbf{1}_{G(M)}(x_1,y_1)\leq\sqrt{\frac
{|M|}{s}}\frac{t}{4x_1^{2-H}}.
\]
Let $\hat{\Psi}_t(x,y)=\sqrt{\frac{|M|}{s}}tx^{H-2}/2$. Then, (\ref
{s3-l4-01}) holds.
Furthermore, by some basic calculations,
%
\begin{eqnarray}\label{s3-l4-4}
&&\int_{G(M)}\phi_t(x_1, y_1)\hat{\Psi}_t(x_1,y_1)\,\mathrm{d}x_1\,\mathrm{d}y_1\qquad\nonumber
\\
&&\quad\leq
\sqrt{\frac{|M|}{s}}\int_{M}^{-t-s}\!\!\!\int_{-t-s-y_1}^\infty\phi_t(x_1,y_1)\frac{t}{2x_1^{2-H}}\,\mathrm{d}x_1\,\mathrm{d}y_1\qquad
\nonumber
\\
&&\quad\leq
\sqrt{\frac{|M|}{s}}\int_{-\infty}^{-t-s}\biggl[\int_{-s-y_1}^\infty\frac{t^2}{2x_1^{4-2H}}\,\mathrm{d}x_1
+\int_{-t-s-y_1}^{-s-y_1}\frac{t(x_1+t+s+y_1)}{2x_1^{4-2H}}\,\mathrm{d}x_1\biggr]\,\mathrm{d} y_1\qquad
\\
&&\quad=
\sqrt{\frac{|M|}{s}}\biggl[\frac{t^{2H}}{4(1-H)(3-2H)}+\int_{0}^{t}\mathrm{d} x_1\int_{-t-s-x_1}^{-t-s}\frac{t(x_1+t+s+y_1)}{2x_1^{4-2H}}\,\mathrm{d}y_1\qquad
\nonumber
\\
&&\quad\hphantom{=\sqrt{\frac{|M|}{s}}\biggl[}{}+
\int_{t}^{\infty}\!\!\mathrm{d}x_1\int_{-t-s-x_1}^{-s-x_1}\frac{t(x_1+t+s+y_1)}{2x_1^{4-2H}}\mathrm{d}y_1\biggr]\qquad
\nonumber
\\
&&\quad=
\frac{1}{4(1-H)(3-2H)(2H-1)}\sqrt{\frac{|M|}{s}}t^{2H}.\qquad\nonumber
\end{eqnarray}
Hence, (\ref{s3-l4-02}) holds for $K_2=\frac
{1}{4(1-H)(3-2H)(2H-1)}$.
\end{pf}

\begin{lem}\label{s3-lem5}
For $0\leq t\leq1$, $M\leq-2-s$, let $\bar{G}(M)=\{-t-s-y<x, y\leq
M\}$ and $\bar{\phi}_t(x,
y)=\phi_t(x,y)\mathbf{1}_{\bar{G}(M)}(x,y)$.
There then exist a non-negative function $\bar{\Psi}_{s, t}(x, y)$ on
$\{x>0, y\leq0\}$
such that for all $n>0$,
%
\begin{equation}\label{s3-l5-01}
F_{n,\bar{\phi}_t}(x_1, y_1)\mathbf{1}_{\bar{G}(M)}(x_1, y_1)\leq\bar
{\Psi}_t(x_1, y_1),
\end{equation}
and a constant $C>0$ which is independent of $t, M$ and $H$, such that
%
\begin{equation}\label{s3-l5-02}
I_2(n,\bar{\phi}_t)\leq\int_{\bar{G}(M)}\phi_t(x_1, y_1)\bar
{\Psi}_t(x_1, y_1)\,\mathrm{d}x_1\,\mathrm{d}y_1\leq K_2Ct^{2H},
\end{equation}
where $K_2$ is the constant in Lemma \ref{s3-lem4}.
\end{lem}

\begin{pf}
From (\ref{s3-la-1}), it follows that
\begin{eqnarray}\label{s3-l5-2}
F_{n,\bar{\phi}_t}(x_1,
y_1)
&=&
n^2\int_{y_1}^{M}\!\!\!\int_{x_1}^{\infty}\bar{\phi}_t(x_2,y_2)\prod_{i=1}^2\sqrt{x_i|y_i|}
\e^{-2n[x_1(y_2-y_1)-(x_2-x_1)y_2]}\,\mathrm{d}x_2\,\mathrm{d}y_2
\nonumber
\\
&=&
\frac{tn^2}{x_1^{3/2-H}}\int_{y_1}^{M}\!\!\sqrt{x_1}|y_1|\int
_{x_1}^{\infty}
\e^{-2n[x_1(y_2-y_1)-(x_2-x_1)y_2]}\,\mathrm{d}x_2\,\mathrm{d}y_2
\nonumber
\\[-10pt]\\[-10pt]
&=&
\frac{t}{2x_1^{2-H}}\int_{|M|}^{|y_1|}\frac{nx_1|y_1|}{|y_2|}
\e^{2nx_1|y_2|}\,\mathrm{d}|y_2|\e^{-2nx_1|y_1|}
\nonumber
\\[-2pt]
&=&
\frac{t}{2x_1^{2-H}}\e^{-2nx_1|y_1|}nx_1|y_1|
\int_{nx_1|M|}^{nx_1|y_1|}\frac{1}{w}\e^{2w}\,\mathrm{d}w.\nonumber
\end{eqnarray}
Let $Q(z)=\mathrm{e}^{-2z}z\int_1^{z}\frac{1}{w}\e^{2w}\,\mathrm{d}w$ for $z\geq1$. Then,
$Q(z)$ is continuous on $[1,\infty)$ and $Q(1)=0$,
$\lim_{z\to\infty}Q(z)=1/2$. Hence, there is a constant $C>0$, which is
independent of $M, t$ and $H$, such that $Q(z)$ is bounded by $C$. Note
that on $\bar{G}(M)$, because $|M|>1$ and $x_1\geq -t-s-M>1,$ we have
$nx_1|M|>1$ for all $n>0$. Therefore, (\ref{s3-l5-2}) yields
\[
F_{n,\bar{\phi}_t}(x_1, y_1)\mathbf{1}_{\bar{G}(M)}(x_1, y_1)\leq\frac
{t}{2x_1^{2-H}}C.
\]
Let $\bar{\Psi}_t(x,
y)=Ctx^{H-2}/2$. By calculations similar to those in (\ref{s3-l4-4}),
(\ref{s3-l5-02})
holds for $K_2=\frac{1}{4(1-H)(3-2H)(2H-1)}$.
\end{pf}

\begin{prop}\label{s3-prop2}
For any $0\leq t\leq1$, there exists a constant $K$,
independent of $t$ but dependent on $s$, such that
%
\begin{equation}\label{s3-p1-0}
\bfE\biggl[\biggl(n\int_0^n\!\!\!\int_{-n}^0
\phi_t(x,y)\sqrt{x|y|}(-1)^{N_{n}(x,y)}\,\mathrm{d}x\,\mathrm{d}y\biggr)^{2}\biggr]
\leq K t^{2H}.
\end{equation}
\end{prop}

\begin{pf}
Using the same notation as in Lemmas
\ref{s3-lem1}--\ref{s3-lem5}, taking $M=-2-s$ and observing that from
(\ref{s3-la-1}), $\phi_t(x, y)=0$ if $(x ,y)$ is not in the set
\[
\{(x, y)\dvtx y<-s, x>0\mbox{ and } x+y>-t-s\}=S\cup G(M)\cup\bar{G}(M),
\]
we can rewrite the left-hand side of (\ref{s3-p1-0}) by
\begin{eqnarray}\label{s3-p1-1}
\hspace*{-30pt}I(n, \phi_t)
:\!&=&
\bfE\biggl[\biggl(n\int_0^n\!\!\!\int_{-n}^0\phi_t(x, y)\sqrt{x|y|}(-1)^{N_{n}(x,y)}\,\mathrm{d}x\,\mathrm{d}y\biggr)^{2}\biggr]\hspace*{35pt}
\nonumber
\\[-10pt]\\[-10pt]
&=&
\bfE\biggl[\biggl(n\int_0^n\!\!\!\int_{-n}^0\bigl(\tilde{\phi }_t(x,y)+\hat{\phi}_t(x, y)+\bar{\phi}_t(x,y)\bigr)
\sqrt{x|y|}(-1)^{N_{n}(x,y)}\,\mathrm{d}x\,\mathrm{d}y \biggr)^{2}\biggr],
\nonumber
\end{eqnarray}
which is bounded by
%
\begin{equation}\label{s3-p1-2}
3I(n, \tilde{\phi}_t)+3I(n,\hat{\phi}_t)+3I(n,\bar{\phi}_t).
\end{equation}
Note that $0\leq \tilde{\phi}_t,\hat{\phi}_t,\bar{\phi}_t\leq\phi_t$.
Lemma \ref{s3-lem1}, Corollary \ref{s3-c1} and Lemma \ref{s3-lema}
imply that
%
\begin{eqnarray}\label{s3-p1-3}
I(n,\tilde{\phi_t})&\leq&2 I_1(n, \tilde{\phi}_t)+2I_2(n,
\tilde{\phi}_t)\nonumber
\\[-10pt]\\[-10pt]
&\leq&\Biggl(\frac {37}{8}+2K_1\sqrt{\frac{1+s}{s}}\Biggr)t^{2H},\nonumber
\end{eqnarray}
and Lemma \ref{s3-lem1}, Corollary \ref{s3-c1} and Lemma \ref{s3-lem4}
imply that
%
\begin{equation}\label{s3-p1-4}
I(n,\hat{\phi}_t)\leq2 I_1(n, \hat{\phi}_t)+2I_2(n,
\hat{\phi}_t)\leq\Biggl(\frac {37}{8}+2K_2\sqrt{\frac{2+s}{s}}\Biggr)t^{2H}.
\end{equation}
Furthermore, from Lemma \ref{s3-lem1}, Corollary \ref{s3-c1} and Lemma
\ref{s3-lem5}, we have
%
\begin{equation}\label{s3-p1-5}
I(n,\bar{\phi}_t)\leq2 I_1(n,
\bar{\phi}_t)+2I_2(n, \bar{\phi}_t)\leq\biggl(\frac{37}{8}+2K_2C\biggr)t^{2H}.
\end{equation}
Therefore, (\ref{s3-p1-1})--(\ref{s3-p1-5}) yield
that
\[
I(n, \phi_t)\leq3\biggl[\frac{111}{8}+2K_1\sqrt{(1+s)/s}+2K_2\sqrt
{(2+s)/s}+2K_2C\biggr]t^{2H}.
\]
Taking $K=3[\frac{111}{8}+2K_1\sqrt{(1+s)/s}+2K_2\sqrt{(2+s)/s}+2K_2C]$
then completes the proof of Proposition \ref{s3-prop2}.
\end{pf}

Finally, we prove Proposition \ref{s3-prop1}.

\begin{pf*}{Proof of Proposition \ref{s3-prop1}}
To prove the tightness of $\{Y_n\}_{n\geq1}$
in ${\mathcal C}[0, 1]$, it suffices to show that for some $r>0$ there
exist two constants $\tilde{M}>0$ and $\eta>1$ such that for any
$t,t'\in[0, 1]$,
\begin{eqnarray*}
&&\bfE\biggl[\biggl(n\int_0^n\!\!\!\int_{-n}^0\frac{g_s(t,x,y)-g_s(t',x,y)}{x^{2-H}}
\sqrt{x|y|}(-1)^{N_{n}(x,y)}\,\mathrm{d}x\,\mathrm{d}y\biggr)^r\biggr]
\\
&&\quad\leq\tilde
{M}(t-t')^{\eta},
\end{eqnarray*}
which follows from the criterion given by Billingsley
(see \cite{b1}, Theorem 12.3) and the fact that our processes are null
at the origin.

Without loss of generality, let $t>t'$. Note that from (\ref{def:g}),
we have
\begin{eqnarray*}
g_s(t,x,y)-g_s(t',x,y)&=&h(t+s, x, y)-h(t'+s, x, y)
\\
&=&g_{t'+s}(t-t', x, y).
\end{eqnarray*}
By Proposition \ref{s3-prop2}, it is easy to check that the above
inequality holds for $r=2$, $\eta=2H>1$ and
$\tilde{M}=3[\frac{111}{8}+2K_1\sqrt{(1+s)/s}+2K_2\sqrt{(2+s)/s)}+2K_2C]$,
which completes the proof of Proposition \ref{s3-prop1}.
\end{pf*}

\begin{rem}\label{rem3.1}
Compared with the proofs of tightness in Bardina et al.~\cite{BJ001,BJ2003}, we have found that under the condition that the kernel $f$ can
be separated by its arguments $(x, y)$, the calculation of $I(n, f)$ is
transformed to the calculation of $I_1(n, f),$ which is relatively
simple; see the proofs of Lemmas 3.1 and 3.2 in \cite{BJ001} and
the proof of Lemma 3.1 in \cite{BJ2003}. However, in our case, the
kernel $f$ cannot be separated by the arguments $(x, y)$, so we need to
discuss $I_2(n, f)$. From our proof, we can see that the calculation of
$I_2(n, f)$ is more complicated and delicate than that of $I_1(n, f)$.
In addition, the fact that the kernel $f$ cannot be separated by the
arguments $(x, y)$ also creates some difficulties in the identification
of the limit law; see the proof of (\ref{s4-p1-8}) in the next section.
\end{rem}

\begin{rem}
By Proposition 10.3 in \cite{EK86}, page 149, Proposition
\ref{s3-prop1} also shows that $\bfP(Y_n\in\mathcal{C}[0,1],
n\in\mathbb{N})=1.$
\end{rem}

\section{Limit law of $Y_n$}\label{4}

In this section, we proceed with the identification of the limit law.
We will prove the following proposition.

\begin{prop}\label{s4-prop1}
The finite-dimensional distributions of $Y_n=\{Y_n(t), t\in[0, 1]\}$
defined by (\ref{pre-3}) converge weakly, as n tends to $\infty$, to
those of a fractional Brownian motion $B^{H}=\{B^{H}(t), t\in[0,1]\}$
with Hurst index $H\in(1/2, 1)$.
\end{prop}

\begin{pf}
For each $k\in\mathbb{N}$,
$a_{1},\ldots,a_{k}\in\mathbb R$ and $0\leq t_1<t_2<\cdots<t_k\leq1$,
we define
\[
L_n=\sum_{j=1}^{k}a_{j}Y_{n}(t_{j})\quad\mbox{and}\quad
U=\sum_{j=1}^{k}a_{j}B^{H}(t_{j}).
\]
It suffices to prove that for any $\xi\in\R$, as $n\to\infty$,
%
\begin{eqnarray}\label{s4-p1-1}
J(n):=|\bfE[\exp(\mathrm{i}\xi L_n)]-\bfE[\exp (\mathrm{i}\xi U)]|\to0.
\end{eqnarray}
For $T>1+s$, let
\[
L_n(T):=n\sum_{j=1}^{k}a_{j}\int_{0}^{T}\!\!\!\int_{-T}^{0} \phi_{t_j}(x,
y)\sqrt{x|y|}(-1)^{N_{n}(x,y)}\,\mathrm{d}x\,\mathrm{d}y
\]
and
\[
U(T):=\sum_{j=1}^{k}a_{j}\int_{0}^{T}\!\!\!\int_{-T}^{0}\phi_{t_j}(x, y)
B(\mathrm{d}x,\mathrm{d}y),
\]
where $B(x,y)$ is given by Lemma \ref{pre-lem-1}. Let
\begin{eqnarray*}
J_1(n,T)
&=&
|\bfE[\exp(\mathrm{i}\xi L_n(T))]-\bfE [\exp(\mathrm{i}\xi U(T))]|,
\\
J_2(n, T)
&=&
|\bfE[\exp(\mathrm{i}\xi L_n)]-\bfE[\exp(\mathrm{i}\xi L_n(T))]|,
\\
J_3(T)
&=&
|\bfE[\exp(\mathrm{i}\xi U(T))]-\bfE[\exp(\mathrm{i}\xi U)]|.
\end{eqnarray*}
Then,
%
\begin{equation}\label{s4-p1-2}
J(n)\leq J_1(n,T)+J_2(n,T)+J_3(T).
\end{equation}
Below, we estimate $J_1(n,T)$, $J_2(n,T)$ and $J_3(T)$, respectively.

(1) We estimate $J_1(n, T)$.

Noting that $\phi_{t_j}(x, y)$ is a non-negative measurable function on
$\{x>0, y\leq0\}$, we can find a sequence of elementary functions
$q^{m,j}(x,y)$ such that
%
\begin{equation}\label{s4-p1-3}
0\leq q^{m,j}(x, y)\leq\phi_{t_j}(x, y)\quad\mbox{and}\quad q^{m, j}(x,
y)\to\phi_{t_j}(x, y)\qquad\mbox{a.e. as }m\to\infty.
\end{equation}
Then, by the dominated convergence theorem, it follows from the fact
that $\int_{0}^{\infty}\!\!\int_{-\infty}^{0}[\phi_{t_j}(x,\break
y)]^2\,\mathrm{d}x\,\mathrm{d}y <\infty$ that as $m\to\infty$,
%
\begin{eqnarray}\label{s4-p1-4}
\int_{0}^{T}\!\!\!\int_{-T}^{0}\bigl(\phi_{t_j}(x,y)-q^{m,j}(x,y)
\bigr)^{2}\,\mathrm{d}x\,\mathrm{d}y\to0.
\end{eqnarray}
For any $m, j, n>0$, define
\begin{eqnarray*}
Y_{n}^{m,j}
&=&
n\int_{0}^{T}\!\!\!\int_{-T}^{0} q^{m,j}(x,y)\sqrt{x|y|}(-1)^{N_{n}(x,y)}\,\mathrm{d}x\,\mathrm{d}y,
\\
B^{m,j}
&=&
\int_{0}^{T}\!\!\!\int_{-T}^{0} q^{m,j}(x,y)B(\mathrm{d}x,\mathrm{d}y).
\end{eqnarray*}
By Lemma \ref{pre-lem-1}, we can readily verify that for fixed $m\in
\mathbb{N}$, as $n\to\infty$,
%
\begin{eqnarray}\label{s4-p1-5}
J_{11}(n, T, m):=\Biggl|\bfE\Biggl[\exp\Biggl(\mathrm{i}\xi
\sum_{j=1}^{k}a_{j}Y_{n}^{m,j}\Biggr)\Biggr]-\bfE\Biggl[\exp\Biggl(\mathrm{i}\xi
\sum_{j=1}^{k}a_{j}B^{m,j}\Biggr)\Biggr]\Biggr|\to0
\end{eqnarray}
because $Y_{n}^{m,j}$ is essentially a linear combination of increments
of $B_n$ defined by (\ref{pre-2}), and $B^{m,j}$ is the same linear
combination of the corresponding limits of increments of $B_n$.

We further define
\begin{eqnarray*}
J_{12}(n, T, m)
&:=&
\Biggl|\bfE[\exp(\mathrm{i}\xi L_n(T))]-\bfE\Biggl[\exp\Biggl(\mathrm{i}\xi\sum_{j=1}^{k}a_{j} Y_{n}^{m,j}\Biggr)\Biggr]\Biggr|,
\\
J_{13}(T,m)
&:=&
\Biggl|\bfE[\exp(\mathrm{i}\xi U(T))]-\bfE\Biggl[\exp\Biggl(\mathrm{i}\xi\sum_{j=1}^{k}a_{j}B^{m,j}\Biggr)\Biggr]\Biggr|.
\end{eqnarray*}
Then, for any $n, m$,
%
\begin{equation}\label{s4-p1-7}
J_1(n, T)\leq J_{11}(n, T, m)+J_{12}(n, T, m)+J_{13}(n, T, m).
\end{equation}
Below, we will show that for all fixed $T$, there exists some
$\gamma_T(m)>0$ such that for any $n>0,$
%
\begin{equation}\label{s4-p1-8}
J_{12}(n, T, m)\leq\gamma_T(m)\to0
\end{equation}
as $m\to\infty$. To this end, let
$f_{m,j}(x,y)=\phi_{t_j}(x,y)-q^{m,j}(x,y)$. Define
\begin{eqnarray*}
\tilde{f}_{m, j}(x, y)&:=&f_{m,j}(x,y)\mathbf{1}_{[0,
T]\times[-t_j-s, 0)}(x, y),
\\ \hat{f}_{m, j}(x, y)&:=&f_{m,j}(x,y)\mathbf{1}_{[0,
T]\times[-T, -t_j-s)}(x, y).
\end{eqnarray*}
By (\ref{s4-p1-3}), as $m\to\infty$,
%
\begin{equation}\label{s4-p1-9}
\phi_{t_j}(x, y)\geq f_{m,j}(x, y)\to0 \qquad\mbox{a.e. in } [0, T]\times[-T,
0].
\end{equation}
Define
\[
R_j(n, m, T)=\bfE\biggl[\biggl|n\int_{0}^{T}\!\!\!\int_{-T}^{0}f_{m,j}(x,
y)\sqrt{x|y|}(-1)^{N_{n}(x,y)}\,\mathrm{d}x\,\mathrm{d}y\biggr|\biggr].
\]
Then, by Lemma
\ref{s3-lem1},
%
\begin{eqnarray}\label{s4-p1-10}
\hspace*{-20pt}[R_j(n, m, T)]^2
&\leq&
n^2E\biggl[\biggl(\int_{0}^{T}\!\!\!\int_{-T}^{0}\bigl(\hat{f}_{m,j}(x,y)+\tilde{f}_{m,j}(x,y)\bigr)\sqrt{x|y|}(-1)^{N_{n}(x,y)}\,\mathrm{d}x\,\mathrm{d}y\biggr)^2\biggr]\quad
\nonumber
\\
&\leq&
2\bigl(I(n, \hat{f}_{m,j})+I(n, \tilde{f}_{m,j})\bigr)
\\
&\leq&
4\bigl(I_1(n, \hat{f}_{m,j})+I_2(n, \hat{f}_{m,j})+I_1(n,\tilde{f}_{m,j})+I_2(n, \tilde{f}_{m,j})\bigr).
\nonumber
\end{eqnarray}
Lemma \ref{s3-lem2} and (\ref{s4-p1-4}) show that for any $n>0$, as
$m\to\infty,$
%
\begin{equation}\label{s4-p1-11}
I_1(n, \tilde{f}_{m,j})\leq\frac{37}{8}\int_0^\infty\!\!\!\int_{-\infty }^0
\tilde{f}_{m, j}^2(x,y)\,\mathrm{d}x\,\mathrm{d}y=:\tilde{\alpha}_j(m,
T)\to0.
\end{equation}
Note that $f_{m, j}(x_1, y_1)\tilde{\Psi}_{t_j}(x_1, y_1)\to0$ a.e.~as
$m\to\infty$. From Lemma \ref{s3-lema} we know that
\[
F_{n,\tilde{f}_{m,j}}(x_1,y_1)\leq F_{n,\tilde{\phi
}_{t_j}}(x_1,y_1)\leq\tilde{\Psi}_{t_j}(x_1, y_1)
\]
and that
\begin{eqnarray*}
&&
\int_0^T\!\!\!\int_{-t_j-s}^0 f_{m, j}(x_1, y_1)F_{n,\tilde{f}_{m,j}}(x_1,y_1)\,\mathrm{d}x_1\,\mathrm{d}y_1
\\
&&\quad\leq
\int_0^T\!\!\!\int_{-t_j-s}^0 \phi_{t_j}(x_1, y_1)\tilde{\Psi}_{t_j}(x_1, y_1)\,\mathrm{d}x_1\,\mathrm{d}y_1<\infty.
\end{eqnarray*}
By the dominated convergence theorem, as $m\to\infty$,
%
\begin{equation}\label{s4-p1-12}
I_2(n, \tilde{f}_{m, j})\leq\tilde{\beta}_j(m,
T):=\int_0^T\!\!\!\int_{-t_j-s}^0 f_{m, j}(x_1, y_1)\tilde{\Psi}_{t_j}(x_1,
y_1)\,\mathrm{d} x_1 \,\mathrm{d}y_1\to0.
\end{equation}
In a similar way, we know there are $\hat{\alpha}_j(m, T)$,
$\hat{\beta}_j(m, T)$ such that for all $n>0$, as $m\to\infty$,
%
\begin{eqnarray}\label{s4-p1-13}
I_1(n,\hat{f}_{m, j})\leq\hat{\alpha}_j(m, T)\to0\quad\mbox{and}\quad
I_2(n,\hat{f}_{m, j})\leq\hat{\beta}_j(m, T)\to0.
\end{eqnarray}
On the other hand, using the mean value theorem, we obtain that
%
\begin{eqnarray}\label{s4-p1-14}
\hspace*{-20pt}J_{12}(n, T, m)
&\leq&
k|\xi| \max_{1\leq j\leq k}\bfE\biggl[\biggl|a_jn\int_{0}^{T}\!\!\!\int_{-T}^0 f_{m,j}(x,y)\sqrt{x|y|}(-1)^{N_{n}(x,y)}\,\mathrm{d}x\,\mathrm{d}y\biggr|\biggr]
\nonumber
\\[-6pt]\\[-6pt]
&=&
k|\xi| \max_{1\leq j\leq k}(|a_j|R_j(n, m, T)).\nonumber
\end{eqnarray}
Hence, (\ref{s4-p1-8}) follows from (\ref{s4-p1-10})--(\ref{s4-p1-14})
with
\[
\gamma_T(m)=2k|\xi|\max_{1\leq j\leq k}\bigl(|a_j|\sqrt{\tilde {\alpha}_j(m,
T) +\tilde{\beta}_j(m, T)+\hat{\alpha}_j(m, T)+\hat{\beta}_j(m, T)}\bigr).
\]

For $J_{13}(m, T)$, we apply the mean value theorem again. Then, as
$m\to\infty$,
%
\begin{eqnarray}\label{s4-p1-15}
J_{13}(m, T)
&\leq&
\xi\bfE\Biggl[\Biggl|a_j\sum_{j=1}^k\int_0^T\!\!\!\int_{-T}^0\bigl(\phi_{t_j}(x,y)-q^{m,j}(x,j)\bigr)B(\mathrm{d}x,\mathrm{d}y)\Biggr|\Biggr]
\nonumber
\\[3pt]
&\leq&
k\xi\max_{1\leq j\leq k}\bfE\biggl[\biggl|a_j\int_0^T\!\!\!\int_{-T}^0\bigl(\phi_{t_j}(x, y)-q^{m,j}(x,j)\bigr)B(\mathrm{d}x,\mathrm{d}y)\biggr|\biggr]
\\[3pt]
&\leq&
k\xi\max_{1\leq j\leq k}\biggl\{\biggl[\int_{0}^{T}\!\!\!\int_{-T}^{0}\bigl(\phi_{t_j}(x, y)-q^{m,j}(x,y)\bigr)^{2}\,\mathrm{d}x\,\mathrm{d}y\biggr]^{1/2}\biggr\}\to0.
\nonumber
\end{eqnarray}

From (\ref{s4-p1-5})--(\ref{s4-p1-8}) and (\ref{s4-p1-15}), we obtain
that for any fixed $T$, as $n\to\infty$,
%
\begin{equation}\label{s4-p1-16}
J_1(n, T)\to0.
\end{equation}

(2) In a similar way as was used to prove (\ref{s4-p1-8}), there
exists some $\zeta(T)>0$ such that for $n>T$, as $T\to\infty$,
%
\begin{equation}\label{s4-p1-17}
J_2(n, T)\leq\zeta(T)\to0.
\end{equation}

(3) In a similar way as was used to prove (\ref{s4-p1-15}), we
obtain that there exists some $\theta(T)>0$ such that as $T\to\infty$,
%
\begin{eqnarray}\label{s4-p1-18}
J_3(T)\leq\theta(T)\to0.
\end{eqnarray}

Therefore, combining (\ref{s4-p1-2}) and
(\ref{s4-p1-16})--(\ref{s4-p1-18}), we can obtain that $J(n)\to0$ as
$n\to\infty$, completing the proof of Proposition \ref{s4-prop1}.
\end{pf}
\eject

Note that Theorem \ref{s4-thhm1} is an immediate result of Propositions
\ref{s3-prop1} and \ref{s4-prop1}. Therefore, the proof of
Theorem \ref{s4-thhm1} is complete.

\section*{Acknowledgments}

This work was done during our visit to the Department of Statistics and
Probability at Michigan State University. The authors thank Professor
Yimin Xiao for stimulating discussions and the department for its good
working conditions. This research was supported in part by the National
Natural Science Foundation of China (No.~10901054). Thanks is also due
to the anonymous referees for their careful reading and detailed
comments which improved the quality of the paper.

\printhistory

\end{document}